\documentclass[11pt]{article}
\usepackage{epsfig,makeidx,amsmath,amsfonts,latexsym}
\newtheorem{theorem}{Theorem}[section]
\newtheorem{conj}[theorem]{Conjecture}

\newtheorem{defn}[theorem]{Definition}

\def\0{{\bf 0}}

\def\eps{\varepsilon}
\def\Cox{{\hfill \Box}}

\begin{document}

\title{The pleasures and pains of studying the two-type Richardson model}

\author{Maria Deijfen  \thanks{Department of Mathematics, Stockholm University,
106 91 Stockholm, Sweden. E-mail: mia@math.su.se} \and Olle H\"{a}ggstr\"{o}m
\thanks{Department of Mathematical Sciences, Chalmers University of
Technology, 412 96 G\"oteborg, Sweden. E-mail: olleh@math.chalmers.se}}

\date{July 2007}

\maketitle

\thispagestyle{empty}

\begin{abstract}

\noindent This paper provides a survey of known results and open
problems for the two-type Richardson model, which is a stochastic
model for competition on $\mathbb{Z}^d$. In its simplest
formulation, the Richardson model describes the evolution of a
single infectious entity on $\mathbb{Z}^d$, but more recently the
dynamics have been extended to comprise two competing growing entities.
For this
version of the model, the main question is whether there is a
positive probability for both entities to simultaneously grow to
occupy infinite parts of the lattice, the conjecture being that
the answer is yes if and only if the entities have the same
intensity. In this paper attention focuses on the two-type model,
but the most important results for the one-type version are also
described.

\vspace{0.5cm}


\vspace{0.5cm}

\end{abstract}

\section{Introduction} 

Consider an interacting particle system in which, at any time $t$,
each site $x\in\mathbb{Z}^d$ is in either of two states, denoted by
0 and 1. A site in state 0 flips to a 1 at rate proportional to the
number of nearest neighbors in state 1, while a site in state 1
remains a 1 forever. We may think of sites in state 1 as being occupied
by some kind of infectious entity, and the model then describes the 
propagation of
an infection where each infected site tries to infect each of its
nearest neighbors on $\mathbb{Z}^d$ at some constant rate $\lambda>0$.
More precisely, if at time $t$ a vertex $x$ is infected and a neighboring
vertex $y$ is uninfected, then, conditional on the dynamics up to time
$t$, the probability that $x$ infects $y$ during a short time
window $(t, t+h)$ is $\lambda h + o(h)$. Here and 
in what follows, sites in state 0 and 1 are referred to as
uninfected and infected respectively. This is the intuitive description 
of the model; a formal definition is given in Section \ref{sect:one-type}. 

The model is a special case of a class of
models introduced by Richardson (1973), and is commonly referred
to as the Richardson model. It has several cousins among processes
from mathematical biology, see e.g.\ Eden (1961), Williams and
Bjerknes (1972) and Bramson and Griffeath (1981). The model is
also a special case of so called first-passage percolation, which
was introduced in Hammersley and Welsh (1965) as a model for
describing the passage of a fluid through a porous medium. In
first-passage percolation, each edge of the $\mathbb{Z}^d$-lattice
is equipped with a random variable representing the time it takes
for the fluid to traverse the edge, and the Richardson model is
obtained by letting these passage times be i.i.d.\ exponential.

Since an infected site stays infected forever, the set of infected
sites in the Richardson model increases to cover all of
$\mathbb{Z}^d$ as $t\to\infty$, and attention focuses on
\emph{how} this set grows. The main result is roughly that the
infection grows linearly in time in each fixed direction and that,
scaled by a factor $1/t$, the set of infected points converges to a
non-random asymptotic shape as $t\to\infty$. To prove that the
growth is linear in a fixed direction involves
Kingman's subadditive ergodic theorem -- in fact, the study of
first-passage percolation was one of the main motivations for the
development of subadditive ergodic theory. That
the linear growth is preserved when all directions are considered
simultaneously is stated in the celebrated shape theorem 
(Theorem \ref{thm:shape} in Section \ref{sect:two-type_bounded})
which originates from Richardson (1973).

Now consider the following extension of the Richardson model,
known as the two-type Richardson model and
introduced in H\"{a}ggstr\"{o}m and Pemantle (1998). Instead of
two possible states for the sites there are three states, which we
denote by 0, 1 and 2. The process then evolves in such a way that, for
$i=1,2$, a site in state 0 flips to state $i$ at rate $\lambda_i$
times the number of nearest neighbors in state $i$ and once in
state 1 or 2, a site remains in that state forever. Interpreting
states 1 and 2 as two different types of infection
and state 0 as absence of infection, this gives rise to a model
describing the simultaneous spread of two infections on
$\mathbb{Z}^d$.  To rigorously define
the model requires a bit more work; see Section \ref{sect:two-type_bounded}.
In what follows we will always assume that $d\geq 2$;
the model makes sense also for $d=1$ but the questions
considered here become trivial.

A number of similar extensions of (one-type) growth models to
(two-type) competition models appear in the literature; see for instance
Neuhauser (1992), Durrett and Neuhauser (1997), Kordzakhia and Lalley (2005)
and Ferrari et al.\ (2006). These tend to require somewhat different 
techniques, and results 
tend not to be easily translated from these other models to the two-type
Richardson model (and vice versa). Closer to the latter are (non-Markovian)
competition models based on first-passage percolation models with 
non-exponential passage time variables -- Garet and Marchand (2005), 
Hoffman (2005:1), Hoffman (2005:2), Garet and Marchand (2006), 
Gou\'er\'e (2007),
Pimentel (2007) -- and a certain continuum model -- Deijfen et al.\ (2004), 
Deijfen and H\"aggstr\"om (2004), Gou\'er\'e (2007). For ease
of exposition, we shall not consider these variations even in cases
where results generalize. 

The behavior of the two-type Richardson model depends on the
initial configuration of the infection and on the ratio between
the intensities $\lambda_1$ and $\lambda_2$ of the infection
types. Assume first, for simplicity, that the model is started at
time 0 from two single sites, the origin being type 1 infected and
the site $(1,0,\ldots,0)$ next to the origin being type 2
infected. Three different scenarios for the development of
the infection are conceivable:

\begin{itemize}
\item[(a)] The type 1 infection at some point completely surrounds
type 2, thereby preventing type 2 from growing any further.

\item[(b)] Type 2 similarly strangles type 1.

\item[(c)] Both infections grow to occupy infinitely many sites.
\end{itemize}

\noindent It is not hard to see that, regardless of the
intensities of the infections, outcomes 
(a) and (b) where one of the infection
types at some point encloses the other have positive
probability regardless of $\lambda_1$ and $\lambda_2$. This is because
each of (a) and (b) can be guaranteed through some finite initial
sequence of infections. In contrast, scenario (c) -- referred to as 
infinite coexistence -- can never be guaranteed from
any finite sequence of infections, and is therefore 
harder to deal with: the main challenge 
is to decide whether, for given values of the parameters
$\lambda_1$ and $\lambda_2$, this event (c) has positive probability
or not. Intuitively, infinite coexistence represents some kind of
power balance between the infections, and it seems reasonable to
suspect that such a balance is possible if and only if the
infections are equally powerful, that is, when
$\lambda_1=\lambda_2$. This is Conjecture \ref{samex_conj} in
Section \ref{sect:two-type_bounded}, which goes back to
H\"{a}ggstr\"{o}m and Pemantle (1998), and, although a lot of
progress have been made, it is not yet fully proved. We describe
the state of the art in Sections \ref{sect:symmetric} and
\ref{sect:nonsymmetric}. 

As mentioned above, apart from the intensities, the development of
the infections in the two-type model also depends on the initial
state of the model. However, if we are only interested in deciding
whether the event of infinite coexistence has positive probability
or not, it turns out that, as long as the initial configuration is
bounded and one of the sets does not completely surround the
other, the precise configuration does not matter, that is, whether
infinite coexistence is possible or not is determined only by the
relation between the intensities. This is proved in Deijfen and
H\"{a}ggstr\"{o}m (2006:1); see Theorem \ref{th:startomr} in
Section \ref{sect:two-type_bounded} for a precise formulation. 
Of course one may also consider unbounded initial configurations.
Starting with both infection types occupying infinitely many sites
means -- apart from in very labored cases -- that they will both
infect infinitely many sites. A more interesting case is when one
of the infection types starts from an infinite set and the other
one from a finite set. We may then ask if outcomes where the
finite type infects infinitely many sites have positive
probability or not. This question is dealt with in Deijfen and
H\"{a}ggstr\"{o}m (2007), and we describe the results in Section
\ref{sect:unbounded}.

The dynamics of the two-type Richardson model is deceptively
simple, and yet gives rise to intriguing 
phenomena on a global scale. In this lies a large part of the pleasure
indicated in the title. Furthermore, 
proofs tend to involve elegant probabilistic
techniques such as coupling, subadditivity and stochastic
comparisons, adding more pleasure. The pain alluded to (which 
by the way is not so severe that it should dissuade readers from
entering this field) comes from the stubborn resistance that some
of the central problems have so far put up against attempts to solve them. 
A case in point is the ``only if'' direction of the aforementioned
Conjecture \ref{samex_conj}, saying that infinite coexistence
starting from a bounded initial configuration does not occur
when $\lambda_1 \neq \lambda_2$. 

\section{The one-type model}  \label{sect:one-type}

As mentioned in the introduction, the one-type Richardson model is
equivalent to first-passage percolation with i.i.d.\ exponential
passage times. To make the construction of the model more precise, first
define $E_{\mathbb{Z}^d}$ as the edge set for the $\mathbb{Z}^d$ lattice
(i.e., each pair of vertices $x,y \in \mathbb{Z}^d$ at Euclidean distance
$1$ from each other have an edge $e \in E_{\mathbb{Z}^d}$ connecting them).
Then attach i.i.d.\ non-negative
random variables $\{\tau(e)\}_{e \in E_{\mathbb{Z}^d}}$ 
to the edges. We take each $\tau(e)$ to be exponentially distributed with
parameter $\lambda>0$, meaning that
\[
P(\tau(e)>t) = \exp(-\lambda t)
\]
for all $t\geq 0$. For $x,y \in \mathbb{Z}^d$, define 
\begin{equation}  \label{eq:path_time}
T(x,y) = \inf_\Gamma \sum_{e \in \Gamma} \tau(e)
\end{equation}
where the infimum is over all paths $\Gamma$ from $x$ to $y$. 
The Richardson model with a given set 
$S_0 \subset \mathbb{Z}^d$ of initially infected sites is now defined 
by taking the set $S_t$ of sites infected at time $t$ to be
\begin{equation}  \label{eq:infected_at_time_t}
S_t = \{ x \in \mathbb{Z}^d : T(y,x) \leq t\mbox{ for some } y \in S_0\} \, .
\end{equation}
It turns out that the infimum in (\ref{eq:path_time}) is a.s.\ a minimum
and attained by a unique path. That $S_t$ grows in the
way described in the introduction is a consequence of the memoryless property
of the exponential distribution: for any $s,t >0$ we have that
$P(\tau(e)>s+t \, | \tau(e)>s) = \exp(-\lambda t)$.

Note that for any $x,y,z \in\mathbb{Z}^d$ we have 
$T(x,y)\leq T(x,z)+T(z,y)$. This subadditivity property
opens up for the use of subadditive ergodic theory in analyzing
the model. To formulate the basic result, let $T(x)$ be the time
when the point $x\in\mathbb{Z}^d$ is infected when starting from a
single infected site at
the origin and write $\mathbf{n}=(n,0,\ldots,0)$. 
It then follows from the subadditive ergodic theorem -- see
e.g.\ Kingman (1968) -- that there is a constant $\mu_\lambda$
such that $T(\mathbf{n})/n\to\mu_\lambda$ almost surely and in
$L_1$ as $n\to\infty$. Furthermore, a simple time scaling argument
implies that $\mu_\lambda=\lambda\mu_1$ and hence, writing
$\mu_1=\mu$, we have that

\begin{equation}\label{eq:time_constant}
\lim_{n\to\infty}\frac{T(\mathbf{n})}{n}=\lambda\mu\quad\textrm{a.s. and in }L_1.
\end{equation}

\noindent The constant $\mu$ indicates the inverse asymptotic
speed of the growth along the axes in a unit rate process and is
commonly referred to as the time constant. It turns out that $\mu>0$,
so that indeed the growth is linear in time. Similarly, 
an analog of (\ref{eq:time_constant}) holds in any direction,
that is, for any $x\in\mathbb{Z}^d$, there is a constant $\mu(x)>0$
such that $T(nx)/n\to\lambda\mu(x)$. The infection hence grows
linearly in time in each fixed direction and the asymptotic speed
of the growth in a given direction is an almost sure constant.

We now turn to the shape theorem, which asserts roughly that the
linear growth of the infection is preserved also when all
directions are considered simultaneously. More precisely, when
scaled down by a factor $1/t$ the set $S_t$ converges to a
non-random shape $A$. To formalize this, let $\tilde{S}_t\subset
\mathbb{R}^d$ be a continuum version of $S_t$ obtained by
replacing each $x\in S_t$ by a unit cube centered at $x$.

\begin{theorem}[Shape Theorem]  \label{thm:shape}
There is a compact convex set $A$ such that, for any $\varepsilon>0$,
almost surely
$$
(1-\varepsilon)\lambda A\subset\frac{\tilde{S}_t}{t} \subset
(1+\varepsilon)\lambda A
$$
for large $t$.
\end{theorem}

\noindent In the above form, the shape theorem was proved in
Kesten (1973) as an improvement
on the original ``in probability" version, which appears
already in Richardson (1973). See also Cox and Durrett (1988) and
Boivin (1991) for generalizations to first-passage percolation
processes with more general passage times. Results concerning
fluctuations around the asymptotic shape can be found, e.g., in Kesten
(1993), Alexander (1993) and Newman and Piza (1995), and, for
certain other passage time distributions, in Benjamini et al.\
(2003).

Working out exactly, or even approximately, what the asymptotic
shape $A$ is has turned out to be difficult. 
Obviously the asymptotic shape inherits all symmetries of the
$\mathbb{Z}^d$ lattice -- invarince under reflection and 
permutation of coordiante hyperplanes -- and it is known to be 
compact and convex,
but, apart from this, not much is known about its qualitative
features. These difficulties with characterizing the shape revolve
around the fact that $\mathbb{Z}^d$ is not rotationally invariant,
which causes the growth to behave differently in different
directions. For instance, simulations on $\mathbb{Z}^2$
indicate that the asymptotic growth is
slightly faster along the axes as compared to the diagonals. There
is however no formal proof of this.

Before proceeding with the two-type model, we mention some work
concerning properties of 
the time-minimizing paths in (\ref{eq:path_time}), also
known as geodesics. Starting at time $0$ with a single infection
at the origin ${\bf 0}$, we denote by $\Gamma(x)$ the (unique) path $\Gamma$
for which the infimum $T({\bf 0}, x)$ in (\ref{eq:path_time}) is
attained. Define
$\Psi=\cup_{x\in\mathbb{Z}^d} \Gamma(x)$, making $\Psi$ a
graph specifying which paths the infection actually takes. It is
not hard to see that $\Psi$ is a tree spanning all of
$\mathbb{Z}^d$ and hence there must be at least one semi-infinite
self-avoiding path from the origin (called an end) in $\Psi$. 
The issue of whether $\Psi$ has more than one end was noted by
H\"aggstr\"om and Pemantle (1998) to
be closely related to the issue of infinite coexistence in the two-type 
Richardson model with $\lambda_1=\lambda_2$: such infinite coexistence
happens with positive probability starting from a finite initial 
configuration if and only if $\Psi$ has at least two ends with
positive probability. 

We say that an infinite path $x_1,x_2,\ldots$ has
asymptotic direction $\hat{x}$ if $x_k/|x_k|\to\hat{x}$ as
$k\to\infty$. In $d=2$, it has been conjectured that every end in
$\Psi$ has an asymptotic direction and that, for every
$x\in\mathbb{R}^2$, there is at least one end (but never more than two) 
in $\Psi$ with asymptotic direction $\hat{x}$. 
In particular, this would mean that $\Psi$ has
uncountably many ends. For results supporting
this conjecture, see Newman (1995) and Newman and Licea (1996).
In the former of these papers, the conjecture is shown to
be true provided an unproven but highly plausible
assumption on the asymptotic shape $A$, saying roughly that the
boundary is sufficiently smooth. See also Lalley (2003) for
related work. 

Results not involving unproven assumptions
are comparatively weak: The coexistence result of H\"aggstr\"om and 
Pemantle (1998) shows for $d=2$ that $\Psi$ has at least two ends with
positive probability. This was later improved to $\Psi$
having almost surely at least $2d$ ends, by
Hoffman (2005:2) for $d=2$ and by Gou\'er\'e (2007) for higher dimensions. 

\section{Introducing two types}
\label{sect:two-type_bounded}

The definition of the two-type Richardson model turns out to be
simplest in the symmetric case $\lambda_1=\lambda_2$, where the
same passage time variables $\{\tau(e)\}_{e \in E_{\mathbb{Z}^d}}$ as 
in the one-type model can be used, with $\lambda= \lambda_1=\lambda_2$. 
Suppose we start with an initial configuration
$(S^1_0, S^2_0)$ of infected sites, where $S^1_0 \subset \mathbb{Z}^d$
are those initially containing type $1$ infection, and 
$S^2_0 \subset \mathbb{Z}^d$
are those initially containing type $2$ infection. We wish to define
the sets $S_t^1$ and $S_t^2$ of type 1 and type $2$ infected sites for
all $t>0$. To this end, set $S_0=S^1_0 \cup
S^2_0$, and take the set $S_t=S_t^1 \cup S_t^2$ of infected sites
at time $t$ to be given by precisely the same formula 
(\ref{eq:infected_at_time_t}) as in the one-type model; a vertex $x\in S_t$
is then assigned infection $1$ or $2$ depending on whether
the $y \in S_0$ for which
\[
\inf\{T(y,x):y\in S_0\}
\]
is attained is in $S_0^1$ or $S_0^2$. 

As in the one-type model, it is a straightforward exercise involving
the memoryless property of the exponential distribution to verify that
$(S_t^1, S_t^2)_{t \geq 0}$ behaves in terms of infection intensities as
described in the introduction. 

This construction demonstrates an intimate link between the one-type
and the symmetric two-type Richardson model: if we watch the two-type model
wearing a pair of of glasses preventing us from distinguishing the two types
of infection, what we see behaves exactly as the one-type model. The link
between infinite coexistence in the two-type model and the number of ends
in the tree of infection $\Psi$ of the one-type model claimed in
the previous section is also a consequence of the construction. 

In the asymmetric case $\lambda_1 \neq \lambda_2$, the two-type model 
is somewhat
less trivial to define due to the fact that the time it takes for
infection to spread along a path depends on the type of infection. 
There are various ways to deal with this, one being to assign, independently
to each $e \in E_{{\mathbb{Z}^d}}$, two independent random variables
$\tau_1(e)$ and $\tau_2(e)$, exponentially distributed with
respective parameters $\lambda_1$ and $\lambda_2$, representing the time
it takes for infections $1$ resp.\ $2$ to traverse $e$. Starting
from an intial configuration $(S^1_0, S^2_0)$, we may picture the infections
as spreading along the edges, taking time $\tau_1(e)$ or $\tau_2(e)$
to cross $e$
depending on the type of infection, with the extra condition that
once a vertex becomes hit by one type of infection it becomes inaccessible
for the other type. This is intuitively clear, but readers with a taste for
detail may require a more rigorous definition, which however we
refrain from here; see H\"aggstr\"om and Pemantle (2000) and
Deijfen and H\"aggstr\"om (2006:1). 

We now move on to describing conjectures
and results. Write $G_i$ for the event that type
$i$ infects infinitely many sites on $\mathbb{Z}^d$ and define
$G=G_1\cap G_2$. The question at issue is:

\begin{equation}\label{eq:coex?}
\textrm{Does $G$ have positive probability?}
\end{equation}

\noindent A priori, the answer to this question may depend both on
the  initial configuration -- that is, on the choice of the sets
$S_0^1$ and $S_0^2$ -- and on the ratio between the
infection intensities $\lambda_1$ and
$\lambda_2$. However, it turns out that, if we
are not interested in the actual value of the probability of $G$,
but only in whether it is positive or not, then the initial
configuration is basically irrelevant, as long as neither of the
initial sets completely surrounds the other. 
This motivates the following definition.\medskip

\begin{defn}
Let $\xi_1$ and $\xi_2$ be two
disjoint finite subsets of $\mathbb{Z}^d$. We say that one of the
sets ($\xi_i$) \emph{strangles} the other ($\xi_j$) if there
exists no infinite self-avoiding path in $\mathbb{Z}^d$ that
starts at a vertex in $\xi_j$ and that does not intersect $\xi_i$.
The pair $(\xi_1,\xi_2)$ is said to be \emph{fertile} if neither
of the sets strangles the other.
\end{defn}

Now write $P^{\lambda_1,\lambda_2}_{\xi_1,\xi_2}$ for
the distribution of a two-type process started from $S_0^1=\xi_1$
and $S_0^2=\xi_2$. We then have the following result.

\begin{theorem}\label{th:startomr}
Let $(\xi_1,\xi_2)$ and $(\xi_1',\xi_2')$ be two fertile pairs of
disjoint finite subsets of $\mathbb{Z}^d$, where $d\geq 2$. For all choices of
$(\lambda_1,\lambda_2)$, we have
$$
P^{\lambda_1,\lambda_2}_{\xi_1,\xi_2}(A)>0\Leftrightarrow
P^{\lambda_1,\lambda_2}_{\xi_1',\xi_2'}(A)>0.
$$
\end{theorem}

For connected initial sets $\xi_1$ and $\xi_2$ and $d=2$, this
result is proved in H\"{a}ggstr\"{o}m and Pemantle (1998). The
idea of the proof in that case is that, by controlling the passage
times of only finitely many edges, two processes started from
$(\xi_1,\xi_2)$ and $(\xi'_1,\xi'_2)$ respectively can be made to
evolve to the same total infected set after some finite time, with
the same configuration of the infection types on the boundary.
Coupling the processes from this time on and observing that the
development of the infections depends only on the boundary
configuration yields the result. This argument however breaks down
when the initial sets are not connected (since it is then not sure
that the same boundary configuration can be obtained in the two
processes) and it is unclear whether it applies for $d\geq 3$. Theorem
\ref{th:startomr} is proved in full generality in Deijfen and
H\"{a}ggstr\"{o}m (2006:1), using a more involved coupling
construction.

It follows from Theorem \ref{th:startomr} that the answer to
(\ref{eq:coex?}) depends only on the value of the intensities
$\lambda_1$ and $\lambda_2$. Hence it is sufficient to consider a
process started from $S_0^1=\mathbf{0}$ and $S_0^2=\mathbf{1}$
(recall that $\mathbf{n}=(n,0,\ldots,0)$),
and in this case we drop subscripts and write $P^{\lambda_1,\lambda_2}$ for
$P^{\lambda_1,\lambda_2}_{{\bf 0}, {\bf 1}}$.
Also, by time-scaling,
we may assume that $\lambda_1=1$. The following conjecture, where
we write $\lambda_2=\lambda$, goes back to H\"{a}ggstr\"{o}m and
Pemantle (1998).

\begin{conj}\label{samex_conj}
In any dimension $d\geq 2$, we have that $P^{1,\lambda}(G)>0$ 
if and only if $\lambda=1$.
\end{conj}

\noindent The conjecture is no doubt true, although
proving it has turned out to be a difficult task. In fact, the
``only if'' direction is not yet fully established. In the
following two sections we describe
the existing results for $\lambda=1$ and $\lambda\neq 1$ respectively.

\section{The case $\lambda=1$}  \label{sect:symmetric}

When $\lambda=1$, we are dealing with two equally powerful
infections and Conjecture \ref{samex_conj} predicts a positive
probability for infinite coexistence. This part of the conjecture 
has been proved:

\begin{theorem}\label{th:lambda=1}
If $\lambda=1$, we have, for any $d\geq 2$, that $P^{1,\lambda}(G)>0$.
\end{theorem}

\noindent This was first proved in the special case $d=2$ by
H\"{a}ggstr\"{o}m and Pemantle (1998). That proof has a very ad hoc flavor,
and heavily exploits not only the two-dimensionality but also
other specific properties of the square lattice,
including a lower
bound on the time constant $\mu$ in (\ref{eq:time_constant})
that just happens to be good enough. When eventually the result
was generalized to higher dimensions, which was done simultaneously
and independently by Garet and Marchand (2005) and Hoffman
(2005:1), much more appealing proofs were obtained. Yet another distinct proof
of Theorem \ref{th:lambda=1} was given by Deijfen and H\"{a}ggstr\"{o}m 
(2007). All four proofs are different, though if you inspect them for
a smallest common denominator you find that they all make critical
use of the fact that the time constant $\mu$ is strictly positive. 
We will give the Garet--Marchand proof below. In Hoffman's proof
ergodic theory is applied to the tree of infection $\Psi$ and a so-called
Busemann function which is shown to exhibit contradictory behavior
under the assumption that infinite coexistence has probability zero. The
Deijfen--H\"{a}ggstr\"{o}m proof proceeds via the two-type Richardson model
with certain infinite initial configurations 
(cf.\ Section \ref{sect:unbounded}). 

\medskip\noindent
{\bf Proof of Theorem \ref{th:lambda=1}:} The following argument
is due to Garet and Marchand (2005), though our presentation follows more
closely the proof of an analogous result in a continuum setting in
Deijfen and H\"{a}ggstr\"{o}m (2004) -- a paper that, despite the 
publication dates, was preceded by and also heavily influenced by
Garet and Marchand (2005). 

Fix a small $\eps>0$. 
By Theorem \ref{th:startomr}, we are free to choose any finite starting
configuration we want, and here it turns out convenient to begin with
a single type $1$ infection at the origin ${\bf 0}$, and a single
type $2$ infection at a vertex ${\bf n}=(n,0,\ldots, 0)$, where
$n$ is large enough so that
\begin{description}
\item{(i)}
$E[T({\bf 0}, {\bf n})] \leq (1 + \eps) n \mu$, and
\item{(ii)}
$P(T({\bf 0}, {\bf n}) < (1 - \eps) n \mu) < \eps$;
\end{description}
note that both (i) and (ii) hold for $n$ large enough due to the
asymptotic speed result (\ref{eq:time_constant}). The reader may
easily check,
for later reference, that (i) and (ii) together with the
nonnegativity of $T({\bf 0}, {\bf n})$ imply for any event $B$
with $P(B)=\alpha$ that
\begin{equation}  \label{eq:key_estimate_GM}
E[T({\bf 0}, {\bf n})\, | \, \neg B] \, 
\leq \, \left(1 + \frac{3\eps}{1-\alpha}\right) n \mu \, .
\end{equation}
Next comes an important telescoping idea: for any positive integer $k$ we have
\begin{eqnarray*}
E[T({\bf 0}, k{\bf n})] & = &
E[T({\bf 0}, {\bf n})]
+ E[T({\bf 0}, 2{\bf n}) - T({\bf 0}, {\bf n})]
+ E[T({\bf 0}, 3{\bf n}) - T({\bf 0}, 2{\bf n})] \\
& & + \ldots +
E[T({\bf 0}, k{\bf n}) - T({\bf 0}, (k-1){\bf n})] \, . 
\end{eqnarray*}
Since $\lim_{k \rightarrow \infty}k^{-1}E[T({\bf 0}, k{\bf n})] = n \mu$,
there must exist arbitrarily large $k$ such that
\[
E[T({\bf 0}, (k+1){\bf n}) - T({\bf 0}, k{\bf n})] \geq (1 - \eps)n \mu \, .
\]
By taking ${\bf m}= k{\bf n}$, and
by translation and reflection invariance, we may deduce that
\begin{equation}  \label{eq:for_arbitrarily_large_m}
E[T({\bf n}, -{\bf m}) - T({\bf 0}, -{\bf m})]  \geq (1 - \eps)n \mu
\end{equation}
for some arbitrarily large $m$. We will pick such an $m$; how large will
soon be specified. 

The goal is to show that $P(G)>0$, so we may assume for contradiction
that $P(G)=0$. By symmetry of the initial configuration, we then have that
$P(G_1)=P(G_2)=\frac{1}{2}$. This implies that
\[
\lim_{m \rightarrow \infty} P({\bf -m} \mbox{ gets infected by type 2})=
\lim_{m \rightarrow \infty} P(T({\bf n}, -{\bf m}) < T({\bf 0}, -{\bf m}))
=\frac{1}{2}
\]
so let us pick $m$ in such a way that
\begin{equation} \label{eq:our_chosen_m}
P(T({\bf n}, -{\bf m}) < T({\bf 0}, -{\bf m})) \geq \frac{1}{4}
\end{equation}
while also (\ref{eq:for_arbitrarily_large_m}) holds. Write $B$ for
the event in (\ref{eq:our_chosen_m}). The expectation
$E[T({\bf n}, -{\bf m}) - T({\bf 0}, -{\bf m})]$ may be
decomposed as
\begin{eqnarray*}
E[T({\bf n}, -{\bf m}) - T({\bf 0}, -{\bf m})] & = &
E[T({\bf n}, -{\bf m}) - T({\bf 0}, -{\bf m}) \, | B]P(B) \\
& & +E[T({\bf n}, -{\bf m}) - T({\bf 0}, -{\bf m}) \, | \neg B]P(\neg B) \\
& \leq & E[T({\bf n}, -{\bf m}) - T({\bf 0}, -{\bf m}) \, | \neg B]P(\neg B) \\
& \leq & \frac{3}{4}
E[T({\bf n}, -{\bf m}) - T({\bf 0}, -{\bf m}) \, | \neg B] \\
& \leq & \frac{3}{4} E[T({\bf n}, {\bf 0}) | \neg B] \\
& \leq & \frac{3}{4} (1 + 4 \eps) n \mu
\end{eqnarray*}
where the second-to-last inequality is due to the triangle
inequality 
$T({\bf n}, -{\bf m}) \leq T({\bf n}, {\bf 0}) + T({\bf 0}, - {\bf m})$,
and the last one uses (\ref{eq:key_estimate_GM}). For small $\eps$,
this contradicts (\ref{eq:for_arbitrarily_large_m}), so the proof is 
complete. $\Cox$

\section{The case $\lambda\neq 1$}  \label{sect:nonsymmetric}

Let us move on to the case when $\lambda\neq 1$, that is, when the
type 2 infection has a different intensity than type 1. It then
seems unlikely that the kind of equilibrium which is necessary for
infinite coexistence to occur would persist in the long run.
However, this part of Conjecture \ref{samex_conj} is not proved.
The best result to date is the following theorem from 
H\"{a}ggstr\"{o}m and Pemantle (2000).

\begin{theorem}\label{th:pain} For any $d\geq 2$, we have
$P^{1,\lambda}(G)=0$ for all but at most countably many values of $\lambda$.
\end{theorem}

\noindent We leave it to the reader to decide whether this
is a very strong or a very weak result: it is very strong in the sense
of showing that infinite coexistence has probability $0$ for
(Lebesgue)-almost all $\lambda$, but very weak in the sense that infinite
coexistence is not ruled out for any given $\lambda$. 

The result may seem a bit peculiar at first sight and
we will spend some time explaining where it comes from and where
the difficulties arise when one tries to strengthen it. Indeed, as
formulated in Conjecture \ref{samex_conj}, the belief is that the
set $\{\lambda:P^{1,\lambda}(G)>0\}$ in fact consists of the
single point $\lambda=1$, but Theorem \ref{th:pain} only asserts
that the set is countable.

First note that, by time-scaling and symmetry, we have 
$P^{1,\lambda}(G)=P^{1,1/\lambda}(G)$ and hence it is enough to
consider $\lambda\leq 1$. An essential ingredient in the proof of
Theorem \ref{th:pain} is a coupling of the two-type processes
$\{P^{1,\lambda}\}_{\lambda\in(0,1]}$ obtained by associating two
independent exponential mean 1 variables $\tau_1(e)$ and
$\tau_2'(e)$ to each edge $e\in\mathbb{Z}^d$ and then letting the
type 2 passage time at parameter value $\lambda$ be given by
$\tau_2(e)=\lambda^{-1}\tau_2'(e)$ and the type 1 time (for any
$\lambda$) by $\tau_1(e)$. Write $Q$ for the probability measure
underlying this coupling and let $G^\lambda$ be the event that
infinite coexistence occurs at parameter value $\lambda$. Theorem
\ref{th:pain} is obtained by showing that
\begin{equation}  \label{eq:in_the_coupling}
\begin{array}{l}
\mbox{with $Q$-probability 1 the event $G^\lambda$ occurs} \\ 
\mbox{for at most one value of $\lambda\in(0,1]$}.
\end{array}
\end{equation}
Hence, $Q(G^\lambda)$ can be positive for at most countably
many $\lambda$, and Theorem
\ref{th:pain} then follows by noting that
$P^{1,\lambda}(G)=Q(G^\lambda)$.

But why is (\ref{eq:in_the_coupling}) true?
Let $G^\lambda_i$ be the event that the type $i$ infection grows
unboundedly at parameter value $\lambda$. Then the coupling
defining $Q$ can be shown to be monotone in the sense that $G^\lambda_1$ is
decreasing in $\lambda$ -- that is, if $G^\lambda_1$ occurs then
$G^{\lambda'}_1$ occurs for all $\lambda'<\lambda$ as well -- and
$G^\lambda_2$ is increasing in $\lambda$. This kind of monotonicity
of the coupling is
crucial for proving (\ref{eq:in_the_coupling}), as is
the following result, which asserts that, on the event that the
type 2 infection survives, the total infected set in a two-type
process with distribution $P^{1,\lambda}$, where $\lambda<1$,
grows to a first approximation
like a one-type process with intensity $\lambda$. More precisely, 
the speed of the growth in the two-type process is determined by
the weaker type 2 infection type. We take $\tilde{S}_t^i$ to denote the
union of all unit cubes centered at points in $S_t^i$ and $A$ is
the limiting shape for a one-type process with rate 1.

\begin{theorem}\label{th:svag_bestammer} Consider a two-type process
with distribution $P^{1,\lambda}$ for some $\lambda\leq 1$. On the
event $G_2$ we have, for any $\varepsilon>0$, that almost surely
$$
(1-\varepsilon)\lambda A\subset\frac{\tilde{S}_t^1\cup
\tilde{S}_t^2}{t}\subset (1+\varepsilon)\lambda A
$$
for large $t$.
\end{theorem}

\noindent Theorem \ref{th:pain} follows readily from this result and the
monotonicity properties of the coupling $Q$. Indeed, fix
$\eps>0$ and suppose $G^\lambda$ occurs. Then Theorem \ref{th:svag_bestammer}
guarantees that on level $\lambda$ the type 1 infection is eventually 
contained in $(1+ \eps)\lambda tA$, a conclusion that
extends to all $\lambda'>\lambda$, because increasing the type 2 infection 
rate does not help type 1. On the other hand, for any
$\lambda'>\lambda$ we get on level $\lambda'$ that
the union of the two infections will -- again by
Theorem \ref{th:svag_bestammer} -- eventually contain
$(1- \eps)\lambda' tA$, so by taking $\eps$ sufficiently small we see that
the type 1 infection is strangled on level $\lambda'$, implying
(\ref{eq:in_the_coupling}), and Theorem \ref{th:pain} follows. 

We will not prove Theorem \ref{th:svag_bestammer}, but 
mention that the hard work in proving it lies in 
establishing a certain key
result (Proposition 2.2 in H\"{a}ggstr\"{o}m and Pemantle (2000))
that asserts that if the strong infection type reaches outside
$(1 + \eps) \lambda tA$ infinitely often, then the weak type is doomed. 
The proof of this uses geometrical arguments, the most important
ingredient being a certain spiral construction, emanating from the part
of the 
strong of infection reaching beyond $(1 + \eps) \lambda tA$, 
and designed to allow the strong type to
completely surround the weak type before the weak type catches up
from inside.

How would one go about to strengthen Theorem \ref{th:pain}
and rule out infinite coexistence for all $\lambda \neq 1$?
One possibility would be to try to derive a
contradiction with Theorem \ref{th:svag_bestammer} from the
assumption that the strong infection type grows unboundedly. For
instance, intuitively it seems likely that the strong type
occupying a positive fraction of the boundary of the infected set
would cause the speed of the growth to exceed the speed prescribed
by the weak infection type. This type of argument is indeed used
in Garet and Marchand (2007) to show, for $d=2$, that on the event
of infinite coexistence the fraction of infected sites occupied
by the strong infection will tend to $0$ as
$t\rightarrow \infty$. This feels like a strong indication that
infinite coexistence does not happen. 

Another approach to strengthening Theorem
\ref{th:pain} in order to obtain the only-if direction of Conjecture
\ref{samex_conj} is based on the
observation that, since coexistence represents a power balance
between the infections, it is reasonable to expect that
$P^{1,\lambda}(G)$ decreases as $\lambda$ moves away from 1. We may
formulate that intuition as a conjecture:

\begin{conj}\label{conj:monotonitet}
For the two-type Richardson model on $\mathbb{Z}^d$ with $d\geq
2$, we have, for $\lambda<\lambda'\in(0,1]$, that
$P^{1,\lambda}(G)\leq P^{1,\lambda'}(G)$.
\end{conj}

\noindent A confirmation of this conjecture would, in combination
with Theorem \ref{th:pain}, clearly establish the only-if
direction of Conjecture \ref{samex_conj}: If $P^{1,\lambda}(G)>0$
for some $\lambda<1$, then, according to Conjecture
\ref{conj:monotonitet}, we would have $P^{1,\lambda'}(G)>0$ for
all $\lambda'\in(\lambda,1]$ as well. But the interval
$(\lambda,1]$ is uncountable, yielding a contradiction to Theorem
\ref{th:pain}. 

Although Conjecture
\ref{conj:monotonitet} might seem close to obvious, it has turned
out to be very difficult to prove. A natural first
attempt would be to use coupling. Consider for instance the
coupling $Q$ described above. As pointed out, the events
$G^\lambda_1$ and $G^\lambda_2$ that the individual infections
grow unboundedly at parameter value $\lambda$ are then monotone in
$\lambda$, but one of them is increasing and the other is decreasing,
so monotonicity of their intersection $G^\lambda$ does not follow. 
Hence more sophisticated arguments are needed.

Observing how our colleagues react during seminars and corridor chat, 
we have noted that it is very tempting to go about trying to prove
Conjecture \ref{conj:monotonitet} by abstract and ``easy'' arguments,
here meaning arguments that do not involve any specifics about the geometry
or graph structure of $\mathbb{Z}^d$. To warn against such attempts,
Deijfen and H\"aggstr\"om (2006:2) constructed graphs on which 
the two-type Richardson model fails to exhibit the monotonicity behavior
predicted in Conjecture \ref{conj:monotonitet}. Let us briefly
explain the results. 

The dynamics of the 
two-type Richardson model can of course be defined on graphs other than 
the $\mathbb{Z}^d$ lattice. 
For a graph $\mathcal{G}$, write Coex$(\mathcal{G})$ for the set
of all $\lambda\geq 1$ such that there exists a finite initial
configuration $(\xi_1, \xi_2)$
for which the two-type Richardson model with infection
intensities $1$ and $\lambda$
started from $(\xi_1, \xi_2)$
yields infinite coexistence with positive
probability. Note that, by time-scaling and interchange of the
infections, coexistence is possible at parameter value $\lambda$
if and only if it is possible at $\lambda^{-1}$, so no
information is lost by restricting to $\lambda\geq 1$. In Deijfen
and H\"aggstr\"om (2006:2) examples of graphs $\mathcal{G}$ are
given that demonstrate that, among others, the following kinds of
coexistence sets Coex$(\mathcal{G})$ are possible:

\begin{itemize}
\item[(i)] Coex$(\mathcal{G})$ may be an interval $(a,b)$ with
$1<a<b$.

\item[(ii)] For any positive integer $k$ the set
Coex$(\mathcal{G})$ may consist of exactly $k$ points.

\item[(iii)] Coex$(\mathcal{G})$ may be countably infinite.

\end{itemize}

\noindent All these phenomena show that the monotonicity suggested
in Conjecture \ref{conj:monotonitet} fails for general graphs.
However, a reasonable guess is that Conjecture
\ref{conj:monotonitet} is true on transitive graphs. Indeed, all
counterexamples provided by Deijfen and H\"aggstr\"om 
are highly non-symmetric (one might even say ugly) 
with certain parts of
the graph being designed specifically with propagation of type 1
in mind, while other parts are meant for type 2. We omit the details. 

\section{Unbounded initial configurations}
\label{sect:unbounded}

Let us now go back to the $\mathbb{Z}^d$ setting and describe some
results from our most recent paper, Deijfen and H\"aggstr\"om (2007), 
concerning the two-type model with unbounded initial configurations. 
Roughly, the model will be started from
configurations where one of the infections occupies a single site
in an infinite ``sea" of the other type. The dynamics is as before
and also the question at issue is the same: can both infection
types simultaneously infect infinitely many sites? With both types
initially occupying infinitely many sites the answer is (apart
from in particularly silly cases) obviously yes, so we will focus
on configurations where type 1 starts with infinitely many sites
and type 2 with finitely many -- for simplicity only one. The
question then becomes whether type 2 is able to survive.

To describe the configurations in more detail, write
$(x_1,\ldots,x_d)$ for the coordinates of a point
$x\in\mathbb{Z}^d$, and define $\mathcal{H}=\{x:x_1=0\}$ and
$\mathcal{L}=\{x:x_1 \leq 0\textrm{ and }x_i=0\textrm{ for }i =2,\ldots, d\}$. 
We will consider the following starting configurations.

\begin{equation}\label{initial_configurations}
\begin{array}{rl}
I(\mathcal{H}): & \mbox{all points in $\mathcal{H}\backslash\{\0\}$ are type 1 infected and}\\
& \0 \mbox{ is type 2 infected, and}\\
I(\mathcal{L}): & \mbox{all points in $\mathcal{L}\backslash\{\0\}$ are type 1 infected and}\\
& \0 \mbox{ is type 2 infected.}
\end{array}
\end{equation}

\noindent Interestingly, it turns out that the set of parameter values for
which type 2 is able to grow indefinitely is slightly different
for these two configuration. First note that, as before, we may
restrict to the case $\lambda_1=1$. Write
$P^{1,\lambda}_{\mathcal{H},\mathbf{0}}$ and
$P^{1,\lambda}_{\mathcal{L},\mathbf{0}}$ for the distribution of
the process started from $I(\mathcal{H})$ and $I(\mathcal{L})$
respectively and with type 2 intensity $\lambda$. The following
result, where $G_2$ denotes the event that type 2 grows
unboundedly, is proved in Deijfen and H\"aggstr\"om (2007).

\begin{theorem}\label{th:unbounded} For the two-type Richardson model in
$d\geq 2$ dimensions, we have

\begin{itemize}
\item[\rm{(a)}] $P^{1,\lambda}_{\mathcal{H},\0}(G_2)>0$ if and only
if $\lambda>1$;
\item[\rm{\rm{(b)}}]$P^{1,\lambda}_{\mathcal{L},\0}(G_2)>0$ if and
only if $\lambda\geq 1$.
\end{itemize}
\end{theorem}

\noindent In words, a strictly stronger type 2 infection will be able to
survive in both configurations, but, when the infections have the
same intensity, type 2 can survive only in the configuration
$I(\mathcal{L})$.

The proof of the if-direction of Theorem \ref{th:unbounded} (a) is
based on a lemma stating roughly that the speed of a hampered
one-type process, living only inside a tube which is bounded in
all directions except one, is close to the speed of an unhampered
process when the tube is large. For a two-type process started
from $I(\mathcal{H})$, this lemma can be used to show that, if the
strong type 2 infection at the origin is successful in the
beginning of the time course, it will take off along the
$x_1$-axis and grow faster than the surrounding type 1
infection inside a tube around the $x_1$-axis, thereby escaping
eradication. The same scenario -- that the type 2 infection rushes
away along the $x_1$-axis -- can, by different means, be proved to
have positive probability in a process with $\lambda=1$ started
from $I(\mathcal{L})$. Infinite growth for type 2 when
$\lambda<1$ is ruled out by the key proposition from H\"aggstr\"om
and Pemantle (2000) mentioned in Section 3. Proving that type 2
cannot survive in a process with $\lambda=1$ started from
$I(\mathcal{H})$ is the most tricky part. The idea is basically to
divide $\mathbb{Z}^d$ in different levels, the $l$-th level being
all sites with $x_1$-coordinate $l$, and then show that the
expected number of type 2 infected sites at level $l$ is constant
and equal to 1. It then follows from a certain comparison with a
one-type process on each level combined with an application of
Levy's 0-1 law that the number of type 2 infected sites at the
$l$-th level converges almost surely to 0 as $l\to\infty$.

Finally we mention a question formulated by Itai Benjamini as well
as by an anonymous referee of Deijfen and H\"aggstr\"om (2007). 
We have seen that, when $\lambda=1$, the type 2 infection at the
origin can grow unboundedly from $I(\mathcal{L})$ but not from
$I(\mathcal{H})$. It is then natural to ask what happens if we
interpolate between these two configurations. More precisely, instead
of letting type 1 occupy only the negative $x_1$-axis (as in
$I(\mathcal{L})$), we let it occupy a cone 
of constant slope around the same
axis. The question then is what the critical slope is for this
cone such that there is a positive probability for type 2 to grow
unboundedly. That type 2 cannot survive when the cone occupies
the whole left half-space follows from Theorem \ref{th:unbounded},
as this situation is equivalent to starting the process from  
$I(\mathcal{H})$. It seems likely, as suggested by Itai Benjamini, 
that this is actually also the critical case, that is, 
infinite growth for type 2 most likely have positive probability 
for any smaller type 1 cone. This however remains to be proved.

\section*{References}

\noindent Alexander, K.\ (1993): A note on some rates of convergence
in first-passage percolation, \emph{Ann. Appl. Probab.} \textbf{3},
81-90.\medskip

\noindent Benjamini, I., Kalai, G. and Schramm, O.\ (2003): First
passage  percolation has sublinear distance variation, \emph{Ann.
Probab.} \textbf{31}, 1970-1978.\medskip

\noindent Bramson, M. and Griffeath, D. (1981): On the
Williams-Bjerknes tumour growth model I, \emph{Ann. Probab.}
\textbf{9}, 173-185.\medskip

\noindent Cox, J.T. and Durrett, R. (1981): Some limit theorems
for percolation processes with necessary and sufficient
conditions, \emph{Ann. Probab.} \textbf{9}, 583-603.\medskip

\noindent Deijfen, M.\ and H\"aggstr\"om, O.\ (2004): Coexistence in a 
two-type continuum growth model, 
\emph{Adv. Appl. Probab.} \textbf{36}, 973-980.\medskip

\noindent Deijfen, M.\ and H\"aggstr\"om, O.\ (2006:1): The initial
configuration is irrelevant for the possibility of mutual unbounded
growth in the two-type Richardson model, {\em Comb. Probab.
Computing} \textbf{15}, 345-353.\medskip

\noindent Deijfen, M.\ and H\"{a}ggstr\"{o}m, O.\ (2006:2):
Nonmonotonic coexistence regions for the two-type Richardson model
on graphs, \emph{Electr. J. Probab.} \textbf{11}, 331-344.\medskip

\noindent Deijfen, M.\ and H\"aggstr\"om, O.\ (2007): The two-type
Richardson model with unbounded initial configurations, \emph{Ann.
Appl. Probab.}, to appear.\medskip

\noindent Deijfen, M., H\"aggstr\"om, O.\  and Bagley, J.\ (2004):
A stochastic model for competing growth on 
$R^d$, \emph{Markov Proc. Relat. Fields} \textbf{10}, 217-248.\medskip

\noindent Durrett, R. (1988): \emph{Lecture Notes on Particle
Systems and Percolation}, Wadsworth $\&$ Brooks/Cole.\medskip

\noindent Durrett, R.\ and Neuhauser, C. (1997): Coexistence results 
for some competition models  \emph{Ann. Appl. Probab.}  \textbf{7}, 
10-45.\medskip

\noindent Eden, M. (1961): A two-dimensional growth process,
\emph{Proceedings of the 4th Berkeley symposium on mathematical
statistics and probability} vol. \textbf{IV}, 223-239, University
of California Press.\medskip

\noindent Ferrari, P., Martin, J.\ and Pimentel, L. (2006), 
Roughening and inclination of competition interfaces, \emph{Phys Rev E}
\textbf{73}, 031602 (4 p). \medskip

\noindent Garet, O.\ and Marchand, R.\ (2005): Coexistence in
two-type first-passage percolation models, {\em Ann. Appl. Probab.}
{\bf 15}, 298-330.\medskip

\noindent Garet, O.\ and Marchand, R.\ (2006): Competition between 
growths governed by Bernoulli percolation,
\emph{Markov Proc. Relat. Fields} \textbf{12}, 695-734.\medskip

\noindent Garet, O.\ and Marchand, R.\ (2007): First-passage competition 
with different speeds: positive density for both species is impossible,
preprint, ArXiV math.PR/0608667.\medskip 

\noindent Gou\'er\'e, J.-B. (2007) Shape of territories in some competing
growth models,  \emph{Ann.
Appl. Probab.}, to appear.\medskip

\noindent H\"{a}ggstr\"{o}m, O.\ and Pemantle, R.\ (1998): First
passage percolation and a model for competing spatial growth,
\emph{J. Appl. Probab.} \textbf{35}, 683-692.\medskip

\noindent H\"{a}ggstr\"{o}m, O.\ and Pemantle, R.\ (2000): Absence
of mutual unbounded growth for almost all parameter values in the
two-type Richardson model, \emph{Stoch. Proc. Appl.} \textbf{90},
207-222.\medskip

\noindent Hammersley, J.\ and Welsh D.\ (1965): First passage
percolation, subadditive processes, stochastic networks and
generalized renewal theory, \emph{1965 Proc. Internat. Res. Semin.,
Statist. Lab., Univ. California, Berkeley}, 61-110,
Springer.\medskip

\noindent Hoffman, C.\ (2005:1): Coexistence for Richardson type
competing spatial growth models, {\em Ann. Appl. Probab.} {\bf 15},
739-747.\medskip

\noindent Hoffman, C.\ (2005:2): Geodesics in first passage percolation, 
preprint, ArXiV math.PR/0508114. \medskip

\noindent Kesten, H.\ (1973): Discussion contribution, \emph{Ann.
Probab.} {\bf 1}, 903. \medskip

\noindent Kesten, H.\ (1993): On the speed of convergence in
first-passage percolation, \emph{Ann. Appl. Probab.} \textbf{3},
296-338.\medskip

\noindent Kingman, J.F.C.\ (1968): The ergodic theory of subadditive 
stochastic processes, \emph{J. Roy. Statist. Soc. Ser. B} \textbf{30}, 
499-510.\medskip

\noindent Kordzakhia, G. and Lalley, S. (2005): A two-species competition 
model on $Z^d$, \emph{Stoch. Proc. Appl.} \textbf{115}, 781-796.\medskip

\noindent Lalley, S. (2003): Strict convexity of the limit shape
in first-passage percolation, {\em Electr. Comm. Probab.} {\bf 8}, 
135--141.\medskip

\noindent Licea, C.\ and Newman, C.\ (1996): Geodesics in
two-dimensional  first-passage percolation, \emph{Ann. Probab.}
\textbf{24}, 399-410.\medskip

\noindent Neuhauser, C.\ (1992): Ergodic theorems for the multitype 
contact process, \emph{Probab. Theory Relat. Fields}  \textbf{91}, 
467-506.\medskip

\noindent Newman, C. (1995): A surface view of first passage
percolation, \emph{Proc. Int. Congr. Mathematicians} \textbf{1,2}
(Zurich 1994), 1017-1023.\medskip

\noindent Newman, C.\ and Piza, M.\ (1995): Divergence of shape
fluctuations in two dimensions, \emph{Ann. Probab.} \textbf{23},
977-1005.\medskip

\noindent Pimentel, L. (2007): Multitype shape theorems for 
first passage percolation models, \emph{Adv. Appl. Probab.}  \textbf{39}, 
53-76.\medskip

\noindent Richardson, D.\ (1973): Random growth in a tessellation,
\emph{Proc. Cambridge Phil. Soc.} \textbf{74}, 515-528.\medskip

\noindent Williams, T. and Bjerknes R. (1972): Stochastic model
for abnormal clone spread through epithelial basal layer,
\emph{Nature} \textbf{236}, 19-21.

\end{document}